%% file: CMelliptic_-_current.tex
\documentclass{amsart}
\input{preamble.tex}

\usepackage[all]{xy}
\usepackage{appendix}
\usepackage{url}
\numberwithin{equation}{section}
\numberwithin{table}{section}

\title[Local constants for CM]{Computing local constants for CM elliptic curves}
\author{Sunil Chetty and Lung Li}
\date{\today}
\thanks{This material is based upon work supported by the National Science Foundation under Grant No. DMS-0757807. A portion of this work was completed as part of the second author's undergraduate capstone research project at Colorado College. The first author would like to thank Colorado College for support during his Riley-scholar post-doctoral fellowship. We would also like to thank the referee for the comments and suggestions, particularly regarding Proposition \ref{newcalc} and Corollary \ref{newcalc-cor}}

\begin{document}
\begin{abstract}
	Let $E/k$ be an elliptic curve with CM by $\Oc$. We determine a formula for (a generalization of) the arithmetic local constant of \cite{MR} at almost all primes of good reduction. We apply this formula to the CM curves defined over $\Q$ and are able to describe extensions $F/\Q$ over which the $\Oc$-rank of $E$ grows.
\end{abstract}
\maketitle
\section{Introduction}
\label{intronotation}
\indent Let $p$ be an odd rational prime. Let $k\subset K\subset L$ be a tower of number fields, with $K/k$ quadratic, $L/K$ $p$-power cyclic, and $L/k$ Galois with a dihedral Galois group, i.e. a lift of $1\neq c\in\Gal{K}{k}$ acts by conjugation on $g\in\Gal{L}{K}$ as $cgc^{-1}=g^{-1}$. In \cite{MR} Mazur and Rubin define arithmetic local constants $\delta_v$, one for each prime $v$ of $K$, which describe the growth in $\Z$-rank of $E$ over the extension $L/K$. Specifically (cf. \cite[Theorem 6.4]{MR}), for $\chi:\Gal{L}{K}\into\Qb^\times$ an injective character and $S$ a set of primes of $K$ containing all primes above $p$, all primes ramified in $L/K$, and all primes where $E$ has bad reduction,
	\begin{equation}
	\label{MR6.4}
		\rk{\Z[\chi]}E(L)^{\chi}-\rk{\Z}E(K)\equiv\sum_{v\in S}\delta_v\mod{2}.
	\end{equation}
To phrase their result this way, we must assume the Shafarevich-Tate Conjecture\footnotemark, and we keep this assumption throughout.
\footnotetext{Without this assumption all statements regarding $\Oc$-rank of $E$ would be replaced by analogous statements regarding $\Oc\otimes\Z_p$-corank of the $p^\infty$-Selmer group $\Selmer{p^\infty}{E}{K}$ of $E$.}\\
\indent In \cite{ChCMAV}, the theory of arithmetic local constants is generalized to address the $\Oc$-rank of varieties with complex multiplication (CM) by an order $\Oc$, and we continue in that direction with specific attention to the elliptic curve case. Following \cite{ChCMAV}, we assume that $\Oc\subset\End{K}{E}$ is the maximal order in a quadratic imaginary field $\bK$, $p$ is unramified in $\Oc$, and $\Oc^c=\Oc^{\dagger}=\Oc$ where $^\dagger$ indicates the action of the Rosati involution (see \cite[\S I.14]{MiAV}). When $\bK\not\subset k$, these assumptions imply $K=k\bK$.\\
\indent Our present aim is to provide a simple formula for the local constants $\delta_v$ (see Definition \ref{constantdefn}) for primes $v\nmid p$ of good reduction. We then will use a result (\cite[\S 6]{ChCMAV}) which generalizes \eqref{MR6.4}, with $\Z$ replaced by $\Oc$, to determine conditions under which the $\Oc$-rank of $E$ grows. In \S\ref{examples} we will describe, via class field theory, dihedral extensions $F/\Q$ which satisfy those conditions, in order to give some concrete setting to the results of \S\ref{computing}.
\section{Computing the local constant}
\label{computing}
\indent Suppose $p$ splits\footnotemark~ in $\Oc$, i.e. $p\Oc=\pp_1\pp_2$, with $\pp_1\neq\pp_2$. We denote $R=\Oc/p\Oc$ and $R_i=\Oc/\pp_i$ for $i=1,2$, so that $R\cong R_1\oplus R_2$.
	\footnotetext{The simpler case of $p$ being inert in $\Oc$, i.e. $\Oc/p\Oc$ is a field, is treated similarly.}
\begin{definition}
\label{rankvector}	
	If $M$ is an $\Oc$-module of exponent $p$, define the $R$-rank of $M$ by
		$$\rk{R}M:=(\rk{R_1}M\otimes_{R}R_1,\rk{R_2}M\otimes_{R}R_2).$$
\end{definition}
\indent The following definition is the same as in \cite{ChCMAV} and \cite{MR}. Fix a prime $v$ of $K$ and let $u$ and $w$ be primes of $k$ below $v$ and of $L$ above $v$, respectively. Denote $k_u$, $K_v$, and $L_w$ for the completions of $k$, $K$, and $L$ at $u$, $v$, and $w$, respectively. If $L_w\neq K_v$, let $L'_w$ be the extension of $K_v$ inside $L_w$ with $[L_w:L'_w]=p$, and otherwise let $L'_w=L_w=K_v$.
\begin{definition}
	\label{constantdefn}	
		Define the arithmetic local constant $\delta_v:=\delta(v,E,L/K)$ by
			$$\delta_v\equiv\rk{R}E(K_v)/(E(K_v)\cap\Norm{L_w}{L'_w}E(L_w))\mod{2}.$$
\end{definition}
Now, we will consider primes $v$ of $K$ such that $E$ has good reduction at $v$, $v\nmid p$, $v^c=v$, and $v$ ramifies in $L/K$ (corresponding to Lemma 6.6 of \cite{MR}). Under these conditions Theorem 5.6 of \cite{MR} shows that
	\begin{equation}
	\label{MR5.6eq}
		\dim{\F_p}E(K_v)/(E(K_v)\cap\Norm{L_w}{L'_w}E(L_w))\equiv\dim{\F_p}E(K_v)[p]\mod{2}.
	\end{equation}
Proposition \ref{MR5.6} below shows that we are able to replace $\dim{\F_p}$ by $\rk{R}$ in \eqref{MR5.6eq}. We first need Lemma \ref{MR5.4-5.5}, which follows Lemmas 5.4-5.5 of \cite{MR}, and our proof is meant only to address the change from $\dim{\F_p}$ to $\rk{R}$.\\
\indent Let $\Kc$ and $\Lc$ be finite extensions of $\Q_\ell$, with $\ell\neq p$, and suppose $\Lc/\Kc$ is a finite extension.
\begin{lemma}
\label{MR5.4-5.5}
	Suppose $\Lc/\Kc$ is cyclic of degree $p$, $E$ is defined over $\Kc$ and has good reduction. Then
		\begin{itemize}
			\item[(i)] $\rk{R}E(\Kc)/pE(\Kc)=\rk{R}E(\Kc)[p]$.
			\item[(ii)] If $\Lc/\Kc$ is ramified then $E(\Kc)/pE(\Kc)=E(\Lc)/pE(\Lc)$ and
				$$\Norm{\Lc}{\Kc}E(\Lc)=pE(\Kc).$$
			\item[(iii)] If $\Lc/\Kc$ is unramified then $\Norm{\Lc}{\Kc}E(\Lc)=E(\Kc).$
		\end{itemize}
\end{lemma}
\begin{proof}
	When $\ell\neq p$ we have $E(\Kc)/pE(\Kc)=E(\Kc)[p^\infty]/pE(\Kc)[p^\infty]$. Since $E(\Kc)[p^\infty]$ is finite, (i) follows from the exact sequence of $\Oc$-modules
		$$0\to E(\Kc)[p]\to E(\Kc)[p^\infty]\to pE(\Kc)[p^\infty]\to 0.$$
The content of (ii) and (iii) is on the level of sets, so the proof is exactly as in Lemma 5.5 of \cite{MR}.
\end{proof}
We return to the notation of Definition \ref{constantdefn}. 
\begin{prop}
\label{MR5.6}
	If $v\nmid p$ and $L_w/K_v$ is nontrivial and totally ramified, then $$\delta_v\equiv\rk{R}E(K_v)[p]\mod{2}.$$
\end{prop}
\begin{proof}
	As in \cite[proof of Thm. 5.6]{MR}, Lemma \ref{MR5.4-5.5}(ii) yields $\Norm{L_w}{L'_w}E(L_w)=pE(L'_w)$ and hence $E(K_v)\cap pE(L'_w)=pE(K_v).$ So by Definition \ref{constantdefn} and Lemma \ref{MR5.4-5.5}(i)
		$$\delta_v\equiv\rk{R}E(K_v)/pE(K_v)\equiv\rk{R}E(K_v)[p]\mod{2}.$$
\end{proof}
Now, fix a prime $v$ of $K$. We denote $\kappa_u$ for the residue field of $k_u$, $q=\#\kappa_u$ for the size of finite field $\kappa_u$, and $\tilde{E}$ for the reduction of $E$ to $\kappa_u$.
\begin{prop}
\label{newcalc}
	Suppose $v\nmid p$, $v$ is ramified in $L/K$, and $v^c=v$. If $E$ has good reduction at $v$, then $\delta_v\equiv (1,1)$ if and only if $p\mid\#\tilde{E}(\kappa_u)$.
\end{prop}
\begin{proof}
	We follow the notation of Lemma 6.6 of \cite{MR}. Since $v^c=v$ we know that $K_v/k_u$ is quadratic, and it is unramified by Lemma 6.5(ii) of \cite{MR}. Let $\Phi$ be the Frobenius generator of $\Gal{K^{ur}_v}{k_u}$, so $\Phi^2$ is the Frobenius of $\Gal{K^{ur}_v}{K_v}$.\\
\indent The proof of Lemma 6.6 of \cite{MR} shows that the product of the eigenvalues $\alpha,\beta$ of $\Phi$ on $E[p]$ is $-1$. Also, $E(K_v)[p]=E[p]^{\Phi^2=1}$ is equal (as a set) to $E[p]$ or is trivial depending on whether or not $\braces{\alpha,\beta}=\braces{1,-1}$, respectively. Since $E$ has CM by $\Oc$, $E[p]$ is a rank 1 $R$-module (see e.g. \cite[\S II.1]{Silv2}), so the former case yields
	$$\delta_v\equiv\rk{R}E(K_v)[p]=(1,1)\mod{2}.$$
\indent By assumption $v\nmid p$, so $p$ is prime to the characteristic of $\kappa_u$, and therefore the reduction map restricted to $p$-torsion is injective (\cite[\S VII.3]{Silv}). We also know $E[p]$ is unramified (\cite[\S VII.4]{Silv}), and so the eigenvalues of $\Phi$ acting on $E[p]$ coincide $\mod{p}$ with the eigenvalues of the $q$-power Frobenius map $\phi_q$ on $\tilde{E}[p]$. We know (\cite[\S V]{Silv}) that the characteristic polynomial of $\phi_q$ is $T^2-aT+q$, where $a=q+1-\#\tilde{E}(\kappa_u)$, and from the above comments $q\equiv -1\mod{p}$. Therefore, $\Phi$ having eigenvalues $\pm 1$ is equivalent to $a\equiv 0\mod{p}$ and in turn equivalent to $\#\tilde{E}(\kappa_u)\equiv 0\mod{p}$.\\
\end{proof}
\begin{cor}
\label{newcalc-cor}	
	If $\bK\not\subset k$ then $\delta_v\equiv (1,1)$.
\end{cor}
\begin{proof}
	To see that $p\mid\#\tilde{E}(\kappa_u)$, we show that $a=0$ under our assumptions on $v$, where $a=q+1-\#\tilde{E}(\kappa_u)$ as above\footnotemark.\footnotetext{That $a=0$ in this case is known (see Exercise 2.30 of \cite{Silv2}, \S 4 Theorem 10 of \cite{LangEF}, or Theorem 7.46 \cite{Shimura} for generalization to higher dimensional abelian varieties), we include an argument for completeness.} The theory of complex multiplication gives $a=\pi_u+\bar{\pi}_u$ for some $\pi_u\in\Oc$ such that $\pi_u\bar{\pi}_u=q$ (see e.g. Theorem 14.16 in \cite{Cox}, \cite[\S II.10]{Silv2}, or \cite{SilverbergGrpOrd} for a thorough discussion). As $\bK\not\subset k$, we have $K=k\bK$, and we let $\psi=\psi_{E/K}$ be the Gr\"{o}ssencharacter associated to $E$ and $K$ (see \cite[\S II.9]{Silv2} or \cite{LangEF}). By comparing their effect on $K(E[\ell])$, where $\ell$ is prime to $v$, we see that $\psi(v)^c=\psi(v^c)$, and since $v=v^c$ we have that $\psi(v)$ is fixed by $c$. It follows that $\psi(v)$ is rational, the corresponding $\pi_v\in\Oc\subset\End{K}{E}$ is integral, and in fact $\pi_v=\pm q$ by degree arguments. In addition $\pi_u^2=\pi_v$, and we will see that $\pi_u=\sqrt{-q}$ is purely imaginary. Indeed, $\pi_u$ having no real part implies $a=\pi_u+\bar{\pi}_u=0$, hence
		$$\#\tilde{E}(\kappa_u)\equiv q+1\equiv 0\mod{p}$$
and $\delta_v\equiv (1,1)$ by Proposition \ref{newcalc}.\\
	\indent Suppose instead that $\pi_u=\sqrt{q}$ is real\footnotemark.\footnotetext{The case $\pi_u=-\sqrt{q}$ follows the same argument.} If in addition we suppose $\pi_u$ is integral then the reduction $\phi_q\in\End{}{\tilde{E}}$ of $\pi_u$ would commute with all endomorphisms of $\tilde{E}$. As $\bK\not\subset k$, there is some $\rho\in\End{K}{E}$ such that $\rho\neq\rho^c$ and hence $\tilde{\rho}\neq\tilde{\rho}^c$. Thus for some $P\in\tilde{E}(\kappa_u)$, $P^c=P$ and $\tilde{\rho}(P^c)\neq\tilde{\rho}^c(P)$. As the action of $c$ on $\kappa_u$ coincides with that of Frobenius $\tilde{\Phi}$, it follows that $\tilde{\rho}$ does not commute with $\tilde{\Phi}$, and in turn $\tilde{\rho}$ does not commute with the Frobenius endomorphism $\phi_q\in\End{}{\tilde{E}}$ induced by $\tilde{\Phi}$.\\
	\indent If $\pi_u=\sqrt{q}$ is real and irrational, then $k\subsetneq\Q(\pi_u)k\subset K$ and so $c\in\Gal{K}{k}$ acts non-trivially on $\pi_u$, i.e. $\pi_u^c=-\sqrt{q}$. It follows that 
		$$q=\Norm{\bK}{\Q}(\pi_u)=\pi_u\pi_u^c=-q,$$
which is a contradiction, and we conclude $\pi_u$ is purely imaginary as desired.
\end{proof}
Define a set $\Sk_L$ of primes $v$ of $K$ by
	$$\Sk_L:=\braces{\text{$v\mid p$, or $v$ ramifies in $L/K$, or where $E$ has bad reduction}}.$$
\begin{theorem}[Theorem 6.1 of \cite{ChCMAV}]
\label{Ch6.1}	
	Let $\chi:\Gal{L}{K}\into\Qb^\times$ be an injective character, and $\Oc[\chi]$ the extension of $\Oc$ by the values of $\chi$. Assuming the Shafarevich-Tate Conjecture,
		$$\rk{\Oc[\chi]}E(L)^{\chi}-\rk{\Oc}E(K)\equiv\sum_{v\in\Sk_L} \delta_v\mod{2}.$$
\end{theorem}
We now consider a dihedral tower $k\subset K\subset F$ where $F/K$ is $p$-power abelian. Following \cite[\S 3]{MR}, we note that there is a bijection between cyclic extensions $L/K$ in $F$ and irreducible rational representations $\rho_L$ of $G=\Gal{F}{K}$. The semisimple group ring $\bK[G]$ decomposes as
	$$\bK[G]\cong\oplus_{L}\bK[G]_L$$
where $\bK[G]_L$ is the $\rho_L$-isotypic component of $\bK[G]$. For each $L$, for us it suffices deal with an injective character $\chi:\Gal{L}{K}\into\Qb^\times$ appearing in the direct-sum decomposition of $\rho_L\otimes\Qb^\times$, and $\rk{\Oc[\chi]}E(F)^{\chi}$ is independent\footnotemark~ of the choice of $\chi$.
	\footnotetext{We could instead write that $\dim{\Qb}(E(F)\otimes\Qb)^\chi$ is independent of the choice of $\chi$.}
\begin{theorem}
\label{mainthm}
	Assume $\bK\not\subset k$.\footnotemark~ Suppose that for every prime $v$ satisfying $v^c=v$ and which ramifies in $F/K$, we have $v\nmid p$ and $E$ has good reduction at $v$. For $m$ equal to the number of such $v$, if $\rk{\Oc}E(K)+m$ is odd then
		$$\rk{\Oc}E(F)\geq [F:K].$$
\end{theorem}
\footnotetext{The case $\bK\subset k$ is similar, with $m$ equal to the number of $v$ such that $p\mid\#\tilde{E}(\kappa_u)$.}
\begin{proof}
	Fix a cyclic extension $L/K$ inside $F$. If $v$ is a prime of $K$ and $v^c\neq v$ then $\delta_v\equiv\delta_{v^c}$ and hence $\delta_v+\delta_{v^c}\equiv (0,0)\mod{2}$ by Lemma 5.1 of \cite{MR}. If $v^c=v$ and $v$ is unramified in $L/K$, then $v$ splits completely in $L/K$ by Lemma 6.5(i) of \cite{MR}. It follows that $\Norm{L_w}{L'_w}$ is surjective and so $\delta_v\equiv (0,0)$ by Definition \ref{constantdefn}. The remaining primes $v$ are precisely those named in the assumption, so Proposition \ref{newcalc-cor} gives $\sum_{v}\delta_v\equiv (m,m)\mod{2}$. Thus,
	$$\rk{\Oc[\chi]}E(L)^\chi\equiv\rk{\Oc}E(K)+m\mod{2}$$
and we have assumed that the right-hand side is odd.\\
\indent From Corollary 3.7 of \cite{MR} it follows that 
	$$\rk{\Oc}E(F)=\sum_L(\dim{\Q}\rho_L)\cdot(\rk{\Oc[\chi]}E(L)^\chi).$$
As the previous paragraph applies for every cyclic $L/K$ in $F$, we see from the decomposition of $\bK[G]$ that $E(F)\otimes\Q$ contains a submodule isomorphic to $\bK[G]$ and the claim follows.
\end{proof}
\section{CM elliptic curves defined over $\Q$}
\label{examples}
\indent Here, we will consider the CM elliptic curves $E$ defined over $\Q$ (as in \cite[A.3]{Silv2}). For each $E$, our aim is to determine\footnotemark\footnotetext{Determined up to the correspondence of class field theory.} examples of dihedral towers $\Q\subset K\subset F$ over which, according to Theorem \ref{mainthm}, the $\Oc$-rank of $E$ grows. As we have assumed $\Oc\subset\End{K}{E}$, we will consider towers in which $K=\bK$ (see \S\ref{intronotation}). All of our calculations will be done using Sage \cite{sage}.\\
%
%
\indent Let $E_D/\Q$ be the elliptic curve of minimal conductor\footnotemark\footnotetext{See p.483 of \cite{Silv2}, with $f=1$ (in Silverman's notation), for a Weierstrauss equation.} defined over $\Q$ with CM by $K_D=\Q(\sqrt{-D})$. We determine computationally\footnotemark\footnotetext{Specifically with Sage's interface to John Cremona's `mwrank' and Denis Simon's `simon\_two\_descent.'} $\rk{\Z}E_D(K_D)$, and for $D=3$ we see that this group is finite. For $D=4,7$, the situation is less certain, as Sage only tells us that $E_D(\Q)$ is finite and $\rk{\Z}E_D(K_D)\leq 2$. For each of the remaining CM curves $E_D$ defined over $\Q$, one can (provably) calculate that $\rk{\Z}E_D(\Q)=1$. We also have that $\rk{\Z}E_D(K_D)\geq\rk{\Z}E_D(\Q)=1$ and $\rk{\Z}E_D(K_D)$ cannot be even, so $\rk{\Oc}E_D(K_D)\geq 1$. For $D=8$, 11, 19, 43, 67, and 163, Sage gives an upper bound\footnotemark[7] of 3 for $\rk{\Z}E_D(K_D)$ and so for these $D$ we can conclude that in fact $\rk{\Oc}E_D(K_D)=1$.
\subsection{Dihedral Extensions of $\Q$}
\indent Recall that $p$ is a fixed odd rational prime. Presently, we also fix $D\in\braces{3,4,7,\ldots,163}$ and let $E=E_D$, $K=K_D$. We are interested in abelian extensions $F/K$ which are dihedral over $\Q$, and these are exactly the extensions contained in the ring class fields of $K$ (see \cite{Cox}, Theorem 9.18).\\
\indent Let $\Oc_f$ be an order in $\Oc_K$ of conductor $f$. We have a simple formula for the class number $h(\Oc_f)$ of $\Oc_f$ using, for example, Theorem 7.24 of \cite{Cox}, and noting that we have $h(\Oc_K)=1$,
		$$h(\Oc_f)=\frac{f}{[\Oc_K^\times:\Oc_f^\times]}
								\cdot\prod_{\text{primes }\ell\mid f}\left(1-\Jac{-D}{\ell}\frac{1}{\ell}\right).$$
For $D\neq 3,4$ we have $\Oc^{\times}_K=\braces{\pm 1}$ and for $D=4$ we have $\#\Oc_K^\times=4$, so in both of these cases $[\Oc_K^\times:\Oc_f^\times]$ is prime to $p$. For $D=3$, one can show that $[\Oc_K^\times:\Oc_f^\times]=3$ when $f>1$. The following paragraphs require only minor adjustments for the case $p=D=3$.\\
\indent Taking $f$ to be an odd rational prime such that $(-D/f)=\pm 1$, the class number becomes $h(\Oc_f)=f\mp 1$ and so the ring class field $H_{\Oc_f}$ associated to $\Oc_f$ is an abelian extension of $K$ of degree $f\mp 1$. Thus, for $f\equiv \pm 1\mod{p}$ we have a (non-trivial) $p$-power subextension $F/K$ which is dihedral over $\Q$.\\
\indent Next we need to understand the ramification in $F/K$.  As $K$ has class number 1, we know there are no unramified extensions of $K$, and so we must ensure that $F$ satisfies the hypotheses of Theorem \ref{mainthm}. A prime $v$ of $K$ ramifies in $H_{\Oc_f}/K$ if and only if $v\mid f\Oc_K$ (see for example exercise 9.20 in \cite{Cox} and recall $f$ is odd). If we choose $f$ so that $-D$ is not a square $\mod{f}$, $f$ is inert in $K/\Q$, and so $f\Oc_K$ is prime and moreover the only prime that ramifies in $H_{\Oc_f}/K$. If $f\Oc_K$ does not ramify in $F/K$ then the $p$-extension $F/K$ is contained in the Hilbert class field $H_K$ of $K$. As $H_K=K$, this is impossible, so $f\Oc_K$ ramifies in $F/K$ and no other primes ramify in $F/K$. Taking $f$ such that $f\nmid D$ and $-D$ is a square $\mod{f}$, we have that $f$ is not inert and does not ramify in $K/\Q$. As in the previous case, the primes of $K$ above $f$ both ramify in the $p$-extension $F/K$ contained in $H_{\Oc_f}$.\\
\indent Now, suppose $\rk{\Oc}E(K)$ is odd\footnotemark\footnotetext{The cases $D=8,11,\ldots,163$ and possibly $D=4,7$.}. To apply Theorem \ref{mainthm}, we must have an even number $m$ of primes $v$ such that $v^c=v$, $v$ ramifies in $F/K$, $E$ has good reducation at $v$ and for which $p\mid\#\tilde{E}(\Z/f\Z)$. First, we can guarantee $m=0$ if the only primes $v$ which ramify in $F/K$ do not satisfy $v^c=v$, e.g. taking $f\nmid D$ with $(-D/f)=1$. Table \ref{tab:nobadprimes} below gives, for each $D$ and for $p=3,5,7$, the smallest prime $f$ which gives an extension of degree $p$ following this recipe. We note that we do not need Proposition \ref{newcalc} for this case.\\
\indent If we wish to allow for primes $v$ satisfying $v^c=v$, we choose two $p$-extensions $F_1$, $F_2$ from two distinct rational primes $f_i$ as above with $f_i\equiv -1\mod{p}$ and $(-D/f_i)=-1$, for $i=1,2$. The compositum $F=F_1F_2$ will satisfy our requirements. Indeed, firstly $F$ is an abelian $p$-extension of $K$ and is contained in the ring class field $H_{\Oc_{f_1f_2}}$, hence dihedral over $\Q$ with only $f_1\Oc_K$ and $f_2\Oc_K$ ramifying in $F/K$. Secondly, as each $f_i$ is inert in $K/\Q$, each is a supersingular prime for $E$ (this follows from the arguments in Corollary \ref{newcalc-cor}) and hence $p$ divides $\#\tilde{E}(\Z/f_i\Z)=f_i+1$. Thus, $E$ and the $p$-extension $F/K$ satisfy the hypotheses of Theorem \ref{mainthm}. Table \ref{tab:twobadprimes} below gives, for each $D$ and for $p=3,5,7$, the smallest pair of distinct primes $f_1, f_2$ which give extensions of degree $p^2$ following this recipe.\\
\indent Next, suppose $\rk{\Oc}E(K)$ is even.\footnotemark\footnotetext{The case $D=3$ and possibly $D=4,7$.} In this case, we need $m$ to be odd in order to apply Theorem \ref{mainthm}. The same ideas as above still work, and in Table \ref{tab:onebadprime} we list, for each $D$ and for $p=3,5,7$, the smallest prime $f$ for which Theorem \ref{mainthm} guarantees rank $\geq p$.
\begin{remark}
	Though there are algorithms in the literature to compute the defining polynomial of a class field (e.g. \cite[\S 6]{Cohen2}, \cite[\S\S 11-3]{Cox}) and such computational problems are of interest independently, we make no attempt here to explicitly determine the ring class fields $H_{\Oc_f}$. As is apparent from Table \ref{tab:twobadprimes}, our method of determining a field to which Theorem \ref{mainthm} applies involves ring class fields of large degree in a computationally inefficient way. 
\end{remark}
\begin{table}
%
	\begin{tabular}{|c||c|c||c|c||c|c|}
		\hline
		 & \multicolumn{2}{|c||}{$p=3$} & \multicolumn{2}{|c||}{$p=5$} & \multicolumn{2}{|c|}{$p=7$}\\
		\hline
		$D$ & $f$ & $[F:K]$ & $f$ & $[F:K]$ & $f$ & $[F:K]$\\
		\hline
		4 & 13 & 3
			& 41 & 5
			& 29 & 7\\
		7 & 43 & 3
			& 11 & 5
			& 29 & 7\\
		8 & 43 & 3
			& 11 & 5
			& 43 & 7\\
		11 & 31 & 3 
			& 31 & 5 
			& 71 & 7\\
		19 & 7 & 3 
			& 11 & 5 
			& 43 & 7\\
		43 & 13 & 3 
			& 11 & 5 
			& 127 & 7\\
		67 & 103 & 3 
			& 71 & 5 
			& 29 & 7\\
		163 & 43 & 3 
			& 41 & 5 
			& 43 & 7\\
		\hline
	\end{tabular}
		\caption{Case $m=0$}
		\label{tab:nobadprimes}
%
	\begin{tabular}{|c||cc|c||cc|c||cc|c|}
		\hline
		 & \multicolumn{3}{|c||}{$p=3$} & \multicolumn{3}{|c||}{$p=5$} & \multicolumn{3}{|c|}{$p=7$}\\
		\hline
		$D$ & $f_1$ & $f_2$ & $[F:K]$ & $f_1$ & $f_2$ & $[F:K]$ & $f_1$ & $f_2$ & $[F:K]$\\
		\hline
		4 & 11 & 23 & 9 
			& 19 & 59 & 25 
			& 83 & 139 & 49\\
		7 & 5 & 41 & 9 
			& 19 & 59 & 25 
			& 13 & 41 & 49\\
		8 & 5 & 23 & 9 
			& 29 & 79 & 25 
			& 13 & 167 & 49\\
		11 & 2 & 29 & 9
			 & 29 & 79 & 25 
			 & 13 & 41 & 49\\
		19 & 2 & 29 & 9
			 & 29 & 59 & 25
			 & 13 & 41 & 49\\
		43 & 2 & 5 & 9
			 & 19 & 29 & 25
			 & 223 & 349 & 49\\
		67 & 2 & 5 & 9
			 & 79 & 109 & 25 
			 & 13 & 41 & 49\\
		163 & 2 & 5 & 9
				& 19 & 29 & 25
				& 13 & 139 & 49\\
		\hline
	\end{tabular}
		\caption{Case $m=2$}
		\label{tab:twobadprimes}
%
	\begin{tabular}{|c||c|c||c|c||c|c|}
		\hline
		 & \multicolumn{2}{|c||}{$p=3$} & \multicolumn{2}{|c||}{$p=5$} & \multicolumn{2}{|c|}{$p=7$}\\
		\hline
		$D$ & $f$ & $[F:K]$ & $f$ & $[F:K]$ & $f$ & $[F:K]$\\
		\hline
		3 & 17 & 3
			& 29 & 5
			& 41 & 7\\
		4 & 11 & 3
			& 19 & 5
			& 83 & 7\\
		7 & 5 & 3
			& 19 & 5
			& 13 & 7\\
		\hline
	\end{tabular}	
		\caption{Case $m=1$}
		\label{tab:onebadprime}
\end{table}
\bibliographystyle{plain}
\bibliography{research}
\vspace{.02in}
\noindent
	\begin{tabular}{ll}
		Sunil Chetty & Lung Li\\
		schetty@csbsju.edu\hspace{0.25in} & leonli319@yahoo.com\\
		Mathematics Department & \\
		College of St. Benedict and St. John's University & \\
	\end{tabular}
\end{document}

%% file: preamble.tex
\usepackage{amsfonts}
\usepackage{amssymb}
\usepackage{amsmath}
\usepackage{amsthm}
\usepackage{amsmath,amscd}

\usepackage[OT2,T1]{fontenc}
\DeclareSymbolFont{cyrletters}{OT2}{wncyr}{m}{n}
\DeclareMathSymbol{\Sha}{\mathalpha}{cyrletters}{"58}

\newtheorem{theorem}{Theorem}[section]
\newtheorem{lemma}[theorem]{Lemma}
\newtheorem{prop}[theorem]{Proposition}
\newtheorem{cor}[theorem]{Corollary}

\theoremstyle{definition}
\newtheorem{definition}[theorem]{Definition}

\newtheorem{remark}[theorem]{Remark}

\newcommand{\Z}{\mathbb{Z}}

\newcommand{\Q}{\mathbb{Q}}

\newcommand{\Qb}{\bar{\mathbb{Q}}}

\newcommand{\F}{\mathbb{F}}

\newcommand{\bK}{\mathbb{K}}

\newcommand{\braces}[1]{\left\{#1\right\}}

\newcommand{\Sk}{\mathfrak{S}}

\newcommand{\Kc}{\mathcal{K}}
\newcommand{\Lc}{\mathcal{L}}

\newcommand{\Oc}{\mathcal{O}}

\newcommand{\chisub}[1]{{\chi_{\lower2.5pt\hbox{$\scriptstyle #1$}}}}

\newcommand{\beq}{\begin{equation}}
\newcommand{\eeq}{\end{equation}}

\newcommand{\beqa}{\begin{eqnarray}}
\newcommand{\eeqa}{\end{eqnarray}}

\newcommand{\beqaN}{\begin{eqnarray*}}
\newcommand{\eeqaN}{\end{eqnarray*}}

\renewcommand{\phi}{\varphi}

\renewcommand{\mod}[1]{\hspace{.025in}\left(\text{mod } #1\right)}

\renewcommand{\dim}[1]{\text{dim}_{#1}}

\newcommand{\rk}[1]{\text{rank}_{#1}}

\newcommand{\into}{\hookrightarrow}

\newcommand{\Jac}[2]{\left(\frac{{#1}}{{#2}}\right)}
\newcommand{\Norm}[2]{N_{{#1}/{#2}}}

\newcommand{\pp}{\mathfrak{p}}

\newcommand{\End}[2]{\text{End}_{#1}(#2)}

\newcommand{\Gal}[2]{\text{Gal}({#1}/{#2})}

\newcommand{\Selmer}[3]{\text{Sel}_{{#1}}({#2}/{#3})}